\documentclass[12pt]{article}
\usepackage{amsthm,amsfonts,amssymb,amscd}

\begin{document}


\begin{center}
\large \bf Canonical and log canonical thresholds\\
of Fano complete intersections
\end{center}\vspace{0.5cm}

\centerline{A.V.Pukhlikov}\vspace{0.5cm}

\parshape=1
3cm 10cm \noindent {\small \quad\quad\quad \quad\quad\quad\quad
\quad\quad\quad {\bf }\newline It is proved that the global
log canonical threshold of a Zariski general Fano complete
intersection of index 1 and codimension $k$ in
${\mathbb P}^{M+k}$ is equal to one, if $M\geqslant 2k+3$
and the maximum of the degrees of defining equations is at
least 8. This is an essential improvements of the previous
results about log canonical thresholds of Fano complete
intersections. As a corollary we obtain the existence of
K\" ahler-Einstein metrics on generic Fano complete intersections
described above.

Bibliography: 18 titles.} \vspace{1cm}

\noindent Key words: Fano variety, log canonical singularity,
hypertangent divisor, K\" ahler-Einstein metric.\vspace{1cm}

\noindent MSC: 14E05, 14E07, 14J45

\section*{Introduction}

{\bf 0.1. Statement of the main result.} The aim of the present paper
is to show that the global (log) canonical threshold of a general
Fano complete intersection of index 1 is at least (respectively, equal)
to one, except for a sufficiently narrow class of Fano complete intersections,
defined by equations of low degree. More precisely, let
$\underline{d}=(d_1,\dots,d_k)$ be an ordered integral vector, where
$k\geqslant 1$ (the value $k$ is {\it not} fixed) and
$2\leqslant d_1\leqslant\dots \leqslant d_k$, and
$$
M=d_1+\dots +d_k-k\geqslant 3.
$$
To every such vector corresponds a family ${\cal F}(\underline{d})$ of
non-singular Fano complete intersections of codimension $k$ in the
complex projective space ${\mathbb P}^{|\underline{d}|}$, where
$|\underline{d}|=d_1+\dots +d_k=M+k$, which we will in the sequel
for simplicity denote by the symbol ${\mathbb P}$:
$$
{\cal F}(\underline{d})\ni F=Q_1\cap\dots\cap Q_k\subset{\mathbb P},
$$
$\mathop{\rm deg} Q_i=d_i$. Obviously, $F$ is a non-singular Fano variety
of index 1, that is, $\mathop{\rm Pic} F={\mathbb Z}K_F$, where
$K_F=-H_F$ is the canonical class of the variety $F$, and $H_F$ is
the class of its hyperplane section in ${\mathbb P}$. The varieties of
the family ${\cal F}(\underline{d})$ are naturally parametrized
by the coefficients of the polynomials, defining the hypersurfaces
$Q_i$.\vspace{0.1cm}

{\bf Conjecture 0.1.} {\it For a general (in the sense of Zariski topology)
variety $F\in {\cal F}(\underline{d})$ and an arbitrary effective divisor
$D\sim nH_F$ on $F$ the pair $(F,\frac{1}{n}D)$ is canonical, that is, for any
exceptional prime divisor $E$ over $F$ the inequality
$$
\mathop{\rm ord}\nolimits_E D\leqslant n\cdot a(E,F),
$$
holds, where $a(E,F)$ is the discrepancy of $E$ with respect to the model}
$F$.\vspace{0.1cm}

The claim of Conjecture 0.1 is usually stated in the following way: the (global)
canonical threshold of the variety $F$ is at least 1. Recall that the
(global) canonical threshold is defined by the equality
$$
\mathop{\rm ct} (F)=\sup\{ \lambda\in {\mathbb Q}_+ \,|\,
(F,\frac{\lambda}{n}D)\,\,\mbox{is canonical for all}\,\, D\in|nH_F|\,\,\mbox{and all}
\,\,n\geqslant 1\},
$$
and the log canonical threshold, respectively, by the equality
$$
\mathop{\rm lct} (F)=\sup\{ \lambda\in {\mathbb Q}_+ \,|\,
(F,\frac{\lambda}{n}D)\,\,\mbox{is log canonical for all}\,\, D\in|nH_F|\,\,\mbox{and all}
\,\,n\geqslant 1\},
$$
The importance of canonical and log canonical thresholds is connected
with their applications to the complex differential geometry and birational
geometry. Tian, Nadel, Demailly and Koll\' ar showed in
\cite{Tian1987,Nadel1990,DemKol2001}, that the inequality
$$
\mathop{\rm lct} (F)>\frac{M}{M+1}
$$
implies the existence of the K\" ahler-Einstein metric on $F$ (this fact
was shown for arbitrary Fano varieties, not only for complete intersections
in the projective space). Since the property of being canonical is stronger
than that of being log canonical, the claim of Conjecture 0.1 implies the
existence of the K\" ahler-Einstein on a general Fano complete intersection
of index 1. This application alone is sufficient to justify the importance
of Conjecture 0.1. For the applications to birational geometry see
Subsection 0.3.\vspace{0.1cm}

Now let us state the main result of the present paper. Let ${\cal D}$ be the
set of ordered integral vectors $\underline{d}$, such that
$2\leqslant d_1\leqslant\dots\leqslant d_k$, $k\geq 2$. For an integer
$a\geqslant 2$ set:
$$
{\cal D}_{\geqslant a}=\{\underline{d}\,|\, d_k\geqslant a\}.
$$
Recall that $|\underline{d}|=d_1+\dots +d_k$.\vspace{0.1cm}

{\bf Theorem 0.1.} {\it Assume that $\underline{d}\in {\cal D}_{\geqslant 8}$
and $|\underline{d}|\geqslant 3k+3$. Then for a Zariski general variety
$F\in {\cal F}(\underline{d})$ the inequality $\mathop{\rm ct}(F)\geqslant 1$
holds.}\vspace{0.1cm}

{\bf Corollary 0.1.} {\it In the assumptions of Theorem 0.1 the equality
$\mathop{\rm lct}(F)=1$ holds, so that on the variety $F$ there is a
K\" ahler-Einstein metric.}\vspace{0.1cm}

Note that the inequality $\mathop{\rm ct}(F)\geqslant 1$ was shown for a general
variety $F\in {\cal F}(\underline{d})$, $\underline{d}\in {\cal D}_{\geqslant 8}$,
under the assumption that $M\geqslant 4k+1$ (that is,
$|\underline{d}|\geqslant 5k+1$), in \cite{Pukh06b}, and under the assumption
that $M\geqslant 3k+4$ (that is, $|\underline{d}|\geqslant 4k+4$), in
\cite{EcklPukh2016}. For more details about the history of this problem see
Subsection 0.5. One should keep in mind that the smaller are the degrees
$d_i$ of the equations defining $F$ (respectively, the higher is the degree
$d=d_1\cdots d_k=\mathop{\rm deg} F$ with the dimension $M=\mathop{\rm dim} F$
fixed), the harder is to prove the inequality $\mathop{\rm ct}(F)=1$. The case
is similar with proving the birational superrigidity of Fano complete
intersections of index 1 \cite{Pukh13a}: the birational superrigidity remains
an open problem in arbitrary dimension only for three types of
complete intersections,
$$
\underline{d}\in\{ (2,\dots,2,2),\, (2,\dots ,2,3),\, (2,\dots ,2,4)\}.
$$
Note also that the canonicity of the pair $(F,\frac{1}{n} D)$ for
{\it any} divisor $D\sim -nK_F$ is a much stronger fact, than
the canonicity of this pair for a {\it general} divisor $D$ of an
{\it arbitrary} mobile linear system $\Sigma\subset |-nK_F|$, and for that
reason it is harder to prove the inequality $\mathop{\rm ct}(F)\geqslant 1$,
than the birationall rigidity.\vspace{0.3cm}


{\bf 0.2. Regular complete intersections.} We understand the condition
that the variety $F\in {\cal F}(\underline{d})$ is Zariski general in
the sense that at every point $o\in F$ the regularity condition (R), which
we will now state, is satisfied. This condition was used in
\cite{Pukh06b,EcklPukh2016}.\vspace{0.1cm}

Let $F=Q_1\cap\dots\cap Q_k\in {\cal F}(\underline{d})$ and
$o\in F$ be an arbitrary point. Fix a system of affine coordinates
$z_*=(z_1,\dots,z_{M+k})$ on ${\mathbb A}^{M+k}\subset {\mathbb P}$
with the origin at the point $o\in {\mathbb A}^{M+k}$. Let
$f_i(z_*)$ be the (non-homogeneous) polynomial defining the
hypersurface $Q_i$ in the affine chart
${\mathbb A}^{M+k}$, $\mathop{\rm deg} F_i=d_i$. Write down
$$
f_i=q_{i,1}+q_{i,2}+\dots +q_{i,d_i},
$$
where $q_{i,j}(z_*)$ are homogeneous polynomials of degree $j$. On
the set $\{q_{i,j}\,|\, 1\leqslant i\leqslant k, 1\leqslant j\leqslant d_i\}$
we introduce the {\it standard order}, setting:

\begin{itemize}

\item $q_{i,j}$ precedes $q_{a,l}$, if $j<l$,

\item $q_{i,j}$ precedes $q_{a,j}$, if $i<a$.

\end{itemize}
Thus placing the polynomials $q_{i,j}$ in the standard order, we get
a sequence of $d_1+\dots +d_k=M+k$ homogeneous polynomials
\begin{equation}\label{03.11.2016.1}
q_{1,1},\, q_{2,1},\,\dots ,\, q_{k,d_k}
\end{equation}
in $M+k$ variables $z_*$.\vspace{0.1cm}

{\bf Definition 0.1.} (i) The complete intersection $F$ is {\it regular
at the point} $o$, if the linear forms $q_{1,1},\dots,q_{k,1}$ are linearly
independent and for any linear form
$$
h(z_*)\not\in \langle q_{1,1},\dots,q_{k,1} \rangle
$$
the sequence of homogeneous polynomials, which is obtained from
(\ref{03.11.2016.1}) by removing the last two polynomials and adding the
form $h$, is regular in ${\cal O}_{o,{\mathbb P}}$.\vspace{0.1cm}

(ii) The complete intersection $F$ satisfies the {\it condition} (R),
if it is regular at every point $o\in F$.\vspace{0.1cm}

In other words, the regularity at the point $o$ means that, removing
from the sequence (\ref{03.11.2016.1}) the last two polynomials and
adding the form $h$, we obtain $M+k-1$ homogeneous polynomials in
$M+k$ variables, the set of common zeros of which is a finite set of
lines in ${\mathbb A}^{M+k}$, passing through the point
$o=(0,\dots,0)$.\vspace{0.1cm}

{\bf Theorem 0.2.} {\it For every tuple
$\underline{d}\in {\cal D}$ there exists a non-empty Zariski open set
${\cal F}_{\rm reg}(\underline{d})\subset {\cal F}(\underline{d})$,
such that every variety $F\in {\cal F}_{\rm reg}(\underline{d})$
satisfies the condition (R).}\vspace{0.1cm}

Now for $\underline{d}\in {\cal D}_{\geqslant 8}$ the claim of
Theorem 0.1 follows from from Theorem 0.2 and the following claim.
\vspace{0.1cm}

{\bf Theorem 0.3.} {\it Assume that
$\underline{d}\in {\cal D}_{\geqslant 8}$. Then for a variety
$F\in {\cal F}_{\rm reg}(\underline{d})$ the inequality
$\mathop{\rm ct}(F)\geqslant 1$ holds.}\vspace{0.3cm}


{\bf 0.3. The canonical threshold and birational rigidity.} Theorem 0.1
has the following application in birational geometry. For an arbitrary
non-singular primitive Fano variety $X$ (that is,
$\mathop{\rm Pic} X={\mathbb Z}K_X$) of dimension $\mathop{\rm dim} X$
define the {\it mobile canonical threshold} $\mathop{\rm mct} (X)$ as
the supremum of such $\lambda\in{\mathbb Q}_+$, that the pair
$(X,\frac{\lambda}{n}D)$ is canonical for a general divisor $D$ of
an arbitrary mobile linear system $\Sigma\subset |-nK_X|$. The inequality
$\mathop{\rm mct} (X)\geqslant 1$ is ``almost equivalent'' to birational
superrigidity of the variety $X$ (for the definition of birational
rigidity and superrigidity see \cite[Chapter 2]{Pukh_book_13a}).\vspace{0.1cm}

In \cite{Pukh05} the following general fact was shown.\vspace{0.1cm}

{\bf Theorem 0.4.} {\it Assume that primitive Fano varieties
$F_1,\dots,F_K$, $K\geqslant 2$, satisfy the conditions $\mathop{\rm
lct} (F_i)=1$ and $\mathop{\rm mct}(F_i) \geqslant 1$. Then their direct
product
$$
V=F_1\times\dots\times F_K
$$
is a birationally superrigid variety. In particular,\vspace{0.1cm}

{\rm (i)} Every structure of a rationally connected fiber space on
the variety $V$ is given by a projection onto a direct factor.
More precisely, let $\beta\colon V^{\sharp}\to S^{\sharp}$ be a
rationally connected fiber space and $\chi\colon V\dashrightarrow
V^{\sharp}$ a birational map. Then there exists a subset of
indices
$$
I=\{i_1,\dots,i_k\}\subset \{1,\dots,K\}
$$
and a birational map
$$
\alpha\colon F_I=\prod\limits_{i\in I}F_i \dashrightarrow
S^{\sharp},
$$
such that the diagram
$$
\begin{array}{rcccl}
& V &\stackrel{\chi}{\dashrightarrow} & V^{\sharp}&\\
\pi_I\!\! &\downarrow & &\downarrow &\!\!\beta\\
& F_I & \stackrel{\alpha}{\dashrightarrow}& S^{\sharp}&
\end{array}
$$
commutes, that is, $\beta\circ\chi=\alpha\circ\pi_I$, where
$$
\pi_I\colon\prod\limits^K_{i=1}F_i\to \prod\limits_{i\in I}F_i
$$
is the natural projection onto a direct factor.\vspace{0.1cm}

{\rm (ii)} Let $V^{\sharp}$ be a variety with ${\mathbb
Q}$-factorial terminal singularities, satisfying the condition
$$
\mathop{\rm dim}\nolimits_{\mathbb Q}(\mathop{\rm
Pic}V^{\sharp}\otimes{\mathbb Q})\leqslant K,
$$
and $\chi\colon V\dashrightarrow V^{\sharp}$ a birational map.
Then $\chi$ is a (biregular) isomorphism.\vspace{0.1cm}

{\rm (iii)} The groups of birational and biregular self-maps of
the variety $V$ coincide:
$$
\mathop{\rm Bir}V=\mathop{\rm Aut}V.
$$
In particular, the group $\mathop{\rm Bir}V$ is
finite.\vspace{0.1cm}

{\rm (iv)} The variety $V$ admits no structures of a fibration
into rationally connected varieties of dimension strictly smaller
than $\mathop{\rm min}\{\mathop{\rm dim}F_i\}$. In particular, $V$
admits no structures of a conic bundle or a fibration into
rational surfaces.\vspace{0.1cm}

{\rm (v)} The variety $V$ is non-rational.}\vspace{0.1cm}

Since the inequality $\mathop{\rm ct}(F) \geqslant 1$ implies that
$\mathop{\rm mct}(F_i) \geqslant 1$ and
$\mathop{\rm lct}(F_i) =1$, Theorem 0.1 implies that generic complete
intersections $F\in {\cal F}(\underline{d})$ with
$\underline{d}\in {\cal D}_{\geqslant 8}$ for
$|\underline{d}|\geqslant 3k+3$ satisfy the assumptions of
Theorem 0.4.\vspace{0.3cm}


{\bf 0.4. The structure of the paper.} In Sections 1-2 we prove Theorem 0.3.
We reproduce the proof sketched in \cite[Section 3.1]{Pukh06b} in full
detail, somewhat modifying the argument given in \cite{Pukh06b},
adjusting it to a wider class of Fano complete intersections. In
principle, the new argument is potentially applicable to proving
the inequality $\mathop{\rm ct}(F) \geqslant 1$ for complete intersections
$F\in {\cal F}(\underline{d})$ with
$\underline{d}\not\in {\cal D}_{\geqslant 8}$.\vspace{0.1cm}

Our main tool is the technique of hypertangent linear systems.
This is a procedure (described in Section 2), the
``input'' of which is an effective divisor $D\sim nH_F$, such that the
pair $(F,\frac{1}{n}D$ is not canonical (under the assumption that
such pairs exist), and the ``output'' of which is an effective 1-cycle $C$
that has a high multiplicity at some point $o\in F$. More precisely, if
$\underline{d}\in {\cal D}_{\geqslant 8}$, then
$\mathop{\rm mult}\nolimits_o C>\mathop{\rm deg} C$, which is impossible.
This contradiction proves Theorem 0.3.\vspace{0.1cm}

In Section 3 we prove Theorem 0.2.\vspace{0.3cm}


{\bf 0.5. Historical remarks and acknowledgements.} As we pointed out
above, the connection between the existence of K\" ahler-Einstein
metrics and the global loc canonical thresholds was established in
\cite{Tian1987,Nadel1990,DemKol2001}. The special importance of those papers
is in that they connected some concepts of complex differential geometry
with some objects of higher-dimensional birational geometry, which makes
it possible to use the results of birational geometry to prove the
existence of K\" ahler-Einstein metrics. That work was started in
\cite{Ch01} and continued in
\cite{ChPark2002,Pukh06b,ChShr2008,Ch2009,Pukh10b,Ch2009b,
ChParkWon2014,EcklPukh2016}. Every time, a computation or estimate
for the global log canonical threshold, obtained by the methods of
birational geometry (the connectedness principle, inversion of
adjunction, the technique of hypertangent divisors) yielded a proof
of existence of K\" ahler-Einstein metrics for new classes of
varieties. Such results are important by themselves, speaking not
of their applications to birational geometry (Theorem 0.4), that is,
of new classes of birationally rigid varieties.\vspace{0.1cm}

Various technical points, related to the constructions of the
present paper, were discussed by the author in his talks given in
2009-2016 at Steklov Mathematical Institute. The author thanks the
members of Divisions of Algebraic Geometry and of Algebra and
Number Theory for the interest to his work. The author also thanks
his colleagues in the Algebraic Geometry research group at the
University of Liverpool for the creative atmosphere and general
support.


\section{Tangent divisors}

In this section we start the proof of Theorem 0.3. We begin
(Subsection 1.1) with some preparatory work: assuming that the pair
$(F,\frac{1}{n}D)$ is not canonical, we show the existence of a
hyperplane section $\Delta$ of the variety $F$, such that the
multiplicity of the restriction of the divisor $D$ onto $\Delta$
at the point $o$ is strictly higher than $2n$. After that
(Subsection 1.2) using the regularity condition (R), we construct
a subvariety $Y\subset\Delta$ of codimension
$(k+1)$ with a high multiplicity at the point $o$.\vspace{0.3cm}

{\bf 1.1. Inversion of adjunction.} Assume that there exists an effective
divisor $D\sim nH_F$ such that the pair $(F,\frac{1}{n}D)$ is
not canonical, that is, there is an exceptional divisor
$E$ over $F$, satisfying the Noether-Fano inequality
$$
\mathop{\rm ord}\nolimits_ED>na(E,F).
$$
By linearity of this inequality in the divisor $D$ (the integer
$n\in{\mathbb Z}_+$ depends linearly on $D$), we may assume that
$D$ is a prime divisor. Let $B\subset F$ be the centre of the
exceptional divisor $E$. It is well known that the estimate
$$
\mathop{\rm mult}\nolimits_BD>n
$$
holds, whence by for example \cite[Proposition 3.6]{Pukh06b}, we
immediately conclude that $\mathop{\rm dim}B\leqslant k-1$.
Consider a point $o\in B$ of general position. Let
$\sigma\colon F^+\to F$ be its blow up,
$E^+=\sigma^{-1}(o)\cong{\mathbb P}^{M-1}$ the exceptional divisor.
For some hyperplane $\Theta\subset E^+$
the inequality
$$
\mathop{\rm mult}\nolimits_oD+\mathop{\rm mult}\nolimits_{\Theta}D^+>2n
$$
holds, where $D^+$ is the strict transform of the divisor $D$ on $F^+$ (see
\cite[Proposition 2.5]{Pukh06b} or
\cite[Chapter 7, Proposition 2.3]{Pukh_book_13a}).\vspace{0.1cm}

Now let us consider a general hyperplane section $\Delta$ of
the complete intersection $F$, containing the point $o$ and
cutting out the hyperplane $\Theta$ on $E^+$ in the sense that
$\Lambda^+\cap E^+=\Theta$. It is easy to see that the restriction
$D_{\Delta}=D|_{\Delta}=(D\circ\Delta)$ of the divisor $D$ on $\Delta$
satisfies the inequality
\begin{equation}\label{05.11.2016.1}
\mathop{\rm mult}\nolimits_oD_{\Delta}>2n.
\end{equation}
The hyperplane section $\Delta$ can be viewed as a complete intersection
of the type $\underline{d}$ in ${\mathbb P}^{M+k-1}$.\vspace{0.3cm}


{\bf 1.2. Intersection with tangent hyperplanes.} Now assume that
$F$ satisfies the condition (R). In the notations of Subsection 0.2
the system of linear equations
$$
q_{1,1}=\dots=q_{k,1}=0
$$
defines the (embedded) tangent space $T_oF\subset T_o{\mathbb P}$.
Obviously, $E^+={\mathbb P}(T_oF)$. Let $h(z_*)$ be the linear form,
defining the hyperplane that cuts out $\Delta$. In particular,
$$
\{h|_{E^+}=0\}=\Theta
$$
and $h\not\in \langle q_{1,1},\dots,q_{k,1}\rangle$. Let
$$
T_i=\{q_{i,1}|_{\Delta}=0\},
$$
$i=1,\dots,k$, be the tangent hyperplane sections of the variety $\Delta$.
By the condition (R), the inequality $\mathop{\rm dim}\Delta \geqslant 2k+2$ and
the Lefschetz theorem (taking into account that the singularities of
the variety $\Delta$ are at most zero-dimensional and $o\in\Delta$
is a non-singular point), we may conclude that for any $j=1,\dots,k$
$$
T_1\cap\dots\cap T_j=(T_1\circ\dots\circ T_j)
$$
is an irreducible subvariety of codimension $j$ in $\Delta$,
which has multiplicity precisely $2^j$ at the point $o$. We will show that
the effective divisor $D_{\Delta}\sim nH_{\Delta}$
(where $H_{\Delta}$ is the class of a hyperplane section of the
complete intersection $\Delta\subset{\mathbb P}^{M+k-1}$), satisfying the
inequality (\ref{05.11.2016.1}), can not exist. Again by the linearity of the
inequality (\ref{05.11.2016.1}) (we will need no other information
about the divisor $D_{\Delta}$), we assume that $D_{\Delta}$ is a prime
divisor. In particular, the inequality (\ref{05.11.2016.1}) implies that
$D_{\Delta}\neq T_1$ (since $\mathop{\rm mult}_oT_1=2$), so that
the effective cycle $(D_{\Delta}\circ T_1)=Y^*_1$ of the scheme-theoretic
intersection of these divisors is well defined and satisfies the inequality
$$
\mathop{\rm mult}\nolimits_oY^*_1>4n,
$$
and moreover, $Y^*_1\sim nH^2_{\Delta}$; in particular,
$$
\frac{\mathop{\rm mult}_oY^*_1}{\mathop{\rm deg}Y^*_1}>\frac{4}{d},
$$
where $d=\mathop{\rm deg}F=d_1\dots d_k$. In the sequel for
simplicity of notations we write
$$
\frac{\mathop{\rm mult}_o}{\mathop{\rm deg}}
$$
for the ratio of multiplicity at the point $o$ to the degree. Let
$Y_1$ be an irreducible component of the cycle $Y^*_1$ with the maximal
value of $\mathop{\rm mult}\nolimits_o/\mathop{\rm deg}$; in particular,
$$
\frac{\mathop{\rm mult}_o}{\mathop{\rm deg}}Y_1>\frac{4}{d}.
$$
Since by construction $Y_1\subset T_1$ and
$$
\frac{\mathop{\rm mult}_o}{\mathop{\rm deg}}(T_1\cap T_2)=\frac{4}{d},
$$
we conclude that $Y_1\not\subset T_2$ and the effective cycle
$(Y_1\circ T_2)=Y^*_2$ is well defined and satisfies the inequality
$$
\frac{\mathop{\rm mult}_o}{\mathop{\rm deg}}Y^*_2>\frac{8}{d}.
$$
Let $Y_2$ be an irreducible component of the cycle $Y^*_2$ with the maximal
value of $\mathop{\rm mult}_o/\mathop{\rm deg}$.\vspace{0.1cm}

Continuing in the same way, we construct a sequence of irreducible subvarieties
$$
D_{\Delta}=Y_0\supset Y_1\supset\dots\supset Y_k
$$
of codimension $\mathop{\rm codim}(Y_j\subset\Delta)=j+1$, satisfying the
inequality
$$
\frac{\mathop{\rm mult}_o}{\mathop{\rm deg}}Y_j>\frac{2}{d}^{j+1}.
$$
The inequality $M\geqslant 2k+3$ is needed to justify the last step
in this construction: by the Lefschetz theorem,
$T_1\cap\dots\cap T_k=(T_1\circ\dots\circ T_k)$ is an irreducible
subvariety of $\Delta$ of codimension $k$, with the multiplicity $2^k$
at the point $o$ and degree $d$, which makes it possible to form the
effective cycle $Y^*_k=(Y_k\circ T_k)$ of codimension $k+1$.\vspace{0.1cm}

We have shown the following claim.\vspace{0.1cm}

{\bf Proposition 1.1.} {\it Assume that the pair $(F,\frac{1}{n}D)$
is not canonical. Then for some point $o\in F$ and a hyperplane
section $\Delta\ni o$, non-singular at the point $o$, there exists
an irreducible subvariety $Y\subset\Delta$ of codimension $k+1$ in $\Delta$,
satisfying the inequality}
\begin{equation}\label{12.11.2016.1}
\frac{\mathop{\rm mult}_o}{\mathop{\rm deg}}Y>\frac{2}{d}^{k+1}.
\end{equation}

In order to complete the proof of Theorem 0.3, we now need the
technique of hypertangent divisors. It is considered in
the next section.


\section{Hypertangent divisors}

In this section we complete the proof of Theorem 0.3. First
(Subsection 2.1) we construct hypertangent linear systems on the variety
$\Delta$ and study their properties. After that (Subsection 2.2) we
select a sequence of general divisors from the hypertangent systems.
Finally, intersecting the subvariety $Y$ with the hypertangent divisors,
we complete the proof of Theorem 0.3 (Subsection 2.3).\vspace{0.3cm}

{\bf 2.1. Hypertangent linear systems.} For $j\in\{1,\dots,d_i\}$
let
$$
f_{i,j}=q_{i,1}+\dots +q_{i,j}
$$
be the truncated equation of the hypersurface $Q_i$. By the symbol
${\cal P}_{a, M+K}$ we denote the linear space of homogeneous
polynomials of degree $a$ in the coordinates
$z_1,\dots,z_{M+k}$. We use this symbol for $a<0$ as well,
setting in that case ${\cal P}_{a, M+K}=\{0\}$.\vspace{0.1cm}

{\bf Definition 2.1.} The linear system of divisors
$$
\Lambda_a=\left\{\left.\left(\sum^k_{i=1}\sum^{d_i-1}_{j=
1}f_{i,j}s_{a-j}\right)\right|_{\Delta},\,\,s_l\in{\cal P}_{l, M+K}\right\}
$$
is the $a$-th hypertangent linear system on $\Delta$ at the point
$o$.\vspace{0.1cm}

Note that by our convention about the negative degrees
only the polynomials $f_{i,j}$ of degree $j\leqslant a$ are
really used in the construction of the system $\Lambda_a$.\vspace{0.1cm}

Set $\delta=d_k$ and for $a\geqslant 2$ set
$$
r_a=\sharp\{i\,|\, 1\leqslant i\leqslant k,\,\,d_i=a\}\in{\mathbb Z}_+.
$$
Obviously, $k_a=0$ for $a\geqslant\delta+1$. The equality
$d_1+\dots+d_k=M+k$ implies that $\delta\leqslant M$. Obviously,
$$
k=k_2+\dots+k_{\delta}.
$$
We say that we are in
\begin{itemize}
       \item the case I, if $k_{\delta}\geqslant 3$,
       \item the case IIA, if $k_{\delta}=2$,
       \item the case IIB, if $k_{\delta}=1$ and $k_{\delta-1}\geqslant 1$,
       \item the case III, if $k_{\delta}=1$ and $k_{\delta-1}=0$.
\end{itemize}
Obviously, one of these cases takes place: we simply listed all
options.\vspace{0.1cm}

For $a\geqslant 2$ set
$$
m_a=\sum_{i\geqslant a}k_i.
$$
It is easy to see that $m_a$ is the number of polynomials of degree $a$
in the seqeunce (\ref{03.11.2016.1}). In the next proposition we sum up
the properties of hypertangent systems that we will need. The symbol
$\mathop{\rm codim}_o$ stands for the codimension in a neighborhood
of the point $o$ with respect to $\Delta$.\vspace{0.1cm}

{\bf Proposition 2.1.} (i) {\it The following inclusion holds:}
$\Lambda_a\subset|aH_{\Delta}|$, {\it where
$H_{\Delta}$ is the class of a hyperplane section of $\Delta$.\vspace{0.1cm}

{\rm (ii)} The following equality holds:
$\mathop{\rm mult}\nolimits_o\Lambda_a=a+1$.\vspace{0.1cm}

{\rm (iii)} In the cases I and IIA for $a=1,\dots,\delta-2$, and in the cases
IIB and III for $a=1,\dots,\delta-3$ the following equality holds:
$$
\mathop{\rm codim}\nolimits_o\mathop{\rm Bs}\Lambda_a=\sum^{a+1}_{i=2}
m_i.
$$

{\rm (iv)} In the case I for $a=\delta-1$, in the cases IIA and IIB for
$a=\delta-2$, and in the case III for $a=\delta-3$ the following equality
holds:} $\mathop{\rm dim} \mathop{\rm Bs}\Lambda_a=1$.\vspace{0.1cm}

Note that the claim (iii) in the case IIA for $a=\delta-2$ and in the case
III for $a=\delta-3$ coincides with the claim (iv) for these cases.\vspace{0.1cm}

{\bf Prooof of Proposition 2.1.} These are the standard facts of the technique
of hypertangent divisors, following immediately from the regularity condition
(Definition 0.1), see \cite[Chapter 3]{Pukh_book_13a}. The claim (i)
is obvious, the claim (ii) follows from the equality
\begin{equation}\label{11.11.2016.1}
f_{i,j}|_{\Delta}=-q_{i,j+1}|_{\Delta}+\dots,
\end{equation}
where the dots stand for the components of degree $j+2$ and higher, and from the
regularity condition. The claims (iii) and (iv) follow from the equality
(\ref{11.11.2016.1}) and the counting of polynomials of degree $j$ in the
sequence (\ref{03.11.2016.1}). For the details, see
\cite[Chapter 3]{Pukh_book_13a}. Q.E.D.\vspace{0.3cm}


{\bf 2.2. Hypertangent divisors.} The next step is constructing a
sequence of hypertangent divisors $D_{i,j}\in\Lambda_i$. From each
hypertangent linear system $\Lambda_i$ we select
$l_i$ divisors, where the integer $l_i$ is defined in the following way:
$l_2=m_3-1$, $l_i=m_{i+1}$ for $i=3,\dots,\delta-3$, finally,

\begin{itemize}

\item in the case I $l_{\delta-2}=m_{\delta-1}$, $l_{\delta-1}=m_{\delta}-2$

\item in the case IIA $l_{\delta-2}=m_{\delta-1}$

\item in the case IIB $l_{\delta-2}=m_{\delta-1}-1$.

\end{itemize}

\noindent For all other values of $i$ set $l_i=0$.\vspace{0.1cm}

Furthermore, for $l_i\neq 0$ we set
$$
{\cal L}_i=\Lambda^{\times l_i}_i
$$
and a tuple of divisors $(D_{i,1},\dots,D_{i,l_i})\in{\cal L}_i$ is denoted
by the symbol $D_{i,*}$. Finally, set
$$
{\cal L}=\prod_{i\geqslant 2}{\cal L}_i,
$$
where the direct product is taken over all $i$ such that $l_i\neq 0$, see
the definition of the integers $l_i$ above. It is easy to see that ${\cal L}$
is the direct product of
$$
\sum_{i\geqslant 2}l_i=M-k-3
$$
factors (precisely the number of polynomials in the sequence
(\ref{03.11.2016.1}), from which all linear and quadratic forms are removed,
together with one cubic polynomial and the last two polynomials).
The elements of of the space ${\cal L}$, that is, the tuples of tuples of divisors
$$
(D_{2,*},D_{3,*},\dots)
$$
are denoted by the symbol $D_{*,*}$.\vspace{0.1cm}

For an arbitrary equidimensional effective cycle $W$ on
$\Delta,\mathop{\rm dim}W\geqslant2$, and a divisor
$D_{i,j}\in\Lambda_i$, such that none of the components of $W$ is contained
in its support $|D_{i,j}|$, we denote by the symbol
$$
(W\circ D_{i,j})_o
$$
the effective cycle of dimension $\mathop{\rm dim}W-1$, which is obtained from
the cycle $(W\circ D_{i,j})$ of the scheme-theoretic intersection of
$W$ and $D_{i,j}$ (see \cite[Chapter 2]{Ful}) by removing all irreducible
components, not containing the point $o$.\vspace{0.3cm}


{\bf 2.3. Proof of Theorem 0.3.} Now everything is ready to apply
the technique of hypertangent systems to the subvariety
$Y\subset\Delta$, constructed in Section 1. The tuple $D_{*,*}\in{\cal L}$
is understood as a tuple of divisors
$$
(D_{2,1},D_{2,2,}\dots,D_{3,1},\dots),
$$
which makes it possible to apply the construction of the
scheme-theoretic intersection at the point $o$, described above,
many times.\vspace{0.1cm}

{\bf Proposition 2.2.} {\it For a general tuple $D_{*,*}\in{\cal L}$
the effective 1-cycle
$$
C=(Y\circ D_{*,*})_o=(\dots((Y\circ D_{2,1})\circ D_{2,2})_o\dots)_o
$$
is well defined and satisfies the inequalities
$$
\mathop{\rm deg}C\leqslant\mathop{\rm deg}Y\cdot\prod_{i\geqslant 2}i^{l_i}
$$
and}
$$
\mathop{\rm mult}\nolimits_oC\geqslant\mathop{\rm mult}\nolimits_oY\cdot
\prod_{i\geqslant 2}(i+1)^{l_i}.
$$

{\bf Proof.} The procedure of constructing the cycle $C$ is justified
by the claims (iii), (iv) of Proposition 2.1, and the inequalities for the
degree and multiplicity follow from the claims (i) and (ii). Q.E.D.\vspace{0.1cm}

Let us prove Theorem 0.3. Assume that $\delta=d_k\geqslant 8$. Combining the
inequality (\ref{12.11.2016.1}) with the inequalities of Proposition
2.2, we obtain the estimate
$$
\frac{\mathop{\rm mult}_o}{\mathop{\rm deg}}C>
\frac{2^{k+1}}{d}\cdot\prod_{i\geqslant 2}\frac{(i+1)^{l_i}}{i^{l_i}},
$$
and after cancellations we see that the inequality
$\mathop{\rm mult}_oC>\mathop{\rm deg}C$ holds. (For the details, see
\cite[Section 3]{Pukh06b}.) This contradiction completes the proof of
Theorem 0.3.


\section{Regular complete intersections}

In this section we prove Theorem 0.2. First (Subsection 3.1), we
reduce the problem to a local problem about violation of the
regularity condition at a fixed point. After that (Subsection 3.2),
we estimate the codimension of the set of tuples of polynomials,
vanishing simultaneously on some line. Finally (Subsection 3.3),
we estimate the codimension of the set of tuples of polynomials,
the set of common zeros of which has an ``incorrect'' dimension,
but is not a line. This completes the proof of Theorem 0.2.\vspace{0.3cm}

{\bf 3.1. Reduction to the local problem.} Following the standard scheme of
proving the regularity conditions (see \cite[Chapter 3]{Pukh_book_13a} or
any paper that makes useof the technique of hypertangent divisors, for
example, \cite{Pukh01} or \cite{EcklPukh2016}), we have to show that
a violation of the local regularity condition at a fixed point $o$
(that is, the condition (i) of Definition 0.1) imposes at least $M+1$
independent conditions on the coefficients of the polynomials
(\ref{03.11.2016.1}). The complete intersection $F$ is non-singular,
so let us fix the linear forms $q_{1,1}\dots,q_{k,1}$ and so the linear
space
$$
T_oF=\{q_{1,1}=\dots =q_{k,1}=0\}.
$$
The last two polynomials in the sequence (\ref{03.11.2016.1}) are not
used in the regularity condition. Let us re-label the polynomials of
the sequence (\ref{03.11.2016.1}), from which all linear forms and the
last two polynomials are removed, by the symbols
$$
p_1,\dots,p_{M-2}.
$$
Now the local regularity condition can be stated as follows:
for any hyperplane $S\subset T_oF$
(in the notations of Definition 0.1 $S=\{h=0\}\cap T_oF$) the sequence
$$
p_1|_S,\dots,p_{M-2}|_S
$$
is regular at the origin. If
${\mathbb S}={\mathbb P}(S)\cong{\mathbb P}^{M-2}$, then this means that
the closed subset
$$
\{p_1|_{\mathbb S}=\dots=p_{M-2}|_{\mathbb S}=0\}
$$
is zero-dimensional. Fix an isomorphism $T_oF\cong{\mathbb C}^M$. Set
$\delta(i)=\mathop{\rm deg}p_i$ and
${\mathbb T}={\mathbb P}(T_oF)\cong{\mathbb P}^{M-1}$. Let
${\cal P}_{a,M}$ be the space of homogeneous polynomials of degree
$a$ on ${\mathbb C}^M$ (or ${\cal P}^{M-1}$) and
$$
{\cal P}_{\mathbb T}=\prod^{M-2}_{i=1}{\cal P}_{\delta(i),M}.
$$
If all polynomials $p_i$ vanish on a line
$L\subset{\mathbb T}$, then, obviously, the local regularity condition
is violated: it is sufficient to take any hyperplane ${\mathbb S}\supset L$.
For that reason the case when the set $\{p_1=\dots=p_{M-2}=0\}$ contains a line
will be considered separately.\vspace{0.3cm}


{\bf 3.2. The case of a line.} Let ${\cal B}^{\rm line}\subset{\cal P}_{\mathbb T}$
be a closed subset of tuples $(p_1,\dots,p_{M-2})$, such that for some
line $L\subset{\mathbb T}$
$$
p_1|_L\equiv\dots\equiv p_{M-2}|_L\equiv 0.
$$

{\bf Proposition 3.1.} {\it The following inequality holds:} $\mathop{\rm codim}
({\cal B}^{\rm line}\subset {\cal P}_{\mathbb T})\geqslant M+1$.\vspace{0.1cm}

{\bf Proof} is obtained by elementary but not quite trivial computations.\vspace{0.1cm}

{\bf Lemma 3.1.} {\it The following inequality holds:}
$$
\mathop{\rm codim}({\cal B}^{\rm line}\subset{\cal P}_{\mathbb T})=
\sum^{M-2}_{i=1}(\delta(i)+1)-2(M-2).
$$

{\bf Proof.} The first component in the right hand side is the codimension
of the set of tuples of polynomials vanishing on a fixed line
$L\subset{\mathbb T}$. Subtracting the dimension of the Grassmanian of lines,
we complete the proof. Q.E.D.\vspace{0.1cm}

Considering the polynomials $q_{i,j}$ for each $i=1,\dots,k$ separately,
we conclude that
$$
\sum^{M-2}_{i=1}\delta(i)=\sum^k_{i=1}\left(\frac{a_i(a_i+1)}{2}-1\right)=
\sum^k_{i=1}\frac{a_i(a_i+1)}{2}-k,
$$
where $a_i=d_i$ for $i=1,\dots,k-2$, $a_{k-1}=a_k=d_k-1$ in the cases I and IIA
and $a_{k-1}=d_{k-1}$, $a_k=d_k-2$ in the cases IIB and III. In any case
$a_i\geqslant 2$ and
$$
a_1+\dots+a_k=M+k-2.
$$

{\bf Lemma 3.2.} {\it The minimum of the quadratic function
$$
\xi(a_1,\dots,a_k)=\sum^k_{i=1}a_i(a_i+1)
$$
on the set of integral vectors $(a_1,\dots,a_k)$ such that all
$a_i\geqslant 2$ and $a_1+\dots+a_k=A$, where $A=ka+l$, $a\in{\mathbb Z}$
and $l\in\{0,1,\dots,k-1\}$, is equal to}
$$
ka^2+(k+2l)a+2l.
$$

{\bf Proof.} Without loss of generality we assume that the set
$(a_1,\dots,a_k)$ is ordered: $a_i\leqslant a_{i+1}$. It is easy to check
that if two positive integers $u,v$ satisfy the inequality
$u\leqslant v-2$, then
$$
u(u+1)+v(v+1)>(u+1)(u+2)+(v-1)v.
$$
Therefore, if $a_i\leqslant a_{i+1}-2$, then, replacing the vector
$\underline{a}=(a_1,\dots,a_k)$ by the vector $\underline{a}'=
(a'_1,\dots,a'_k)$, where $a'_j=a_j$ for $j\neq i,i+1$, $a'_i=a_i+1$
and $a'_{i+1}=a_{i+1}-1$, we decrease the value of the function $\xi$.
Similarly, if
$$
a_i+1=a_{i+1}=\dots=a_{i+\alpha}=a_{i+\alpha+1}-1,
$$
then, replacing the vector $\underline{a}$ by the vector $\underline{a}'$ with
$$
a'_i=a'_{i+1}=\dots=a'_{i+\alpha+1}=a_{i+1},
$$
we decrease the value of the function $\xi$. In both cases the vector
$\underline{a}'$ remains ordered and satisfies the restrictions
$a'_1+\dots+a'_k=A$, $a'_i\geqslant 2$. Since the set of such vectors
is finite, applying finitely many modifications of the two types
described above, we obtain a vector with
$$
a_1=\dots=a_{k-l}=a\quad\mbox{and}\quad a_{k-l+1}=\dots=a_k=a+1,
$$
which realizes the minimum of the function $\xi$. Simple computations
complete the proof of the lemma. Q.E.D.\vspace{0.1cm}

Now writing $M-2=ka+l$ with
$l\in\{0,\dots,k-1\}$ and applying Lemmas 3.1 and 3.2, after obvious
simplifications we obtain the inequality
$$
\mathop{\rm codim}({\cal B}^{\rm line})\subset{\cal P})\geqslant
\frac12(k(a+1)^2+k(a-1)+2l(a+1)-(M-2).
$$
Now the inequality of Proposition 3.1 follows from the estimate
$$
k(a^2-a)+2l(a-1)\geqslant 6,
$$
which is easy to check (recall that by assumption
$M\geqslant 2k+3$, so that $ak+l\geqslant 2k+1$ and, in particular,
$a\geqslant 2$). Q.E.D. for Proposition 3.1.\vspace{0.1cm}

Starting from this moment, we assume that the polynomials
$p_1,\dots,p_{M-2}$ do not vanish simultaneously on a line
$L\subset{\mathbb T}$.\vspace{0.3cm}


{\bf 3.3. End of the proof of Theorem 0.2.} Fix a hyperplane
${\mathbb S}\subset{\mathbb T}$ and its isomorphism
${\mathbb S}\cong{\mathbb P}^{M-2}$. Set
$$
{\cal P}=\prod^{M-2}_{i=1}{\cal P}_{\delta(i),M-1}.
$$
Since the hyperplane ${\mathbb S}$ varies in a $(M-1)$-dimensional
family, it is sufficient to show that the codimension of the set of tuples
$(p_1,\dots,p_{M-2})\in{\cal P}$ such that the closed set
$$
\{p_1=\dots=p_{M-2}=0\}
$$
has a component of positive dimension, {\it which is not a line},
is of codimension at least $(M+1)+(M-1)=2M$ in ${\cal P}$. Let us check this
fact. The check is not difficult, arguments of that type were published
many times in full detail, so we will just sketch the main steps.\vspace{0.1cm}

Let ${\cal B}_i\subset{\cal P}$ be the set of such tuples that the
closed set
\begin{equation}\label{23.11.2016.1}
\{p_1=\dots=p_{i-1}=0\}\subset{\mathbb P}^{M-2}
\end{equation}
(if $i=0$, then this set is assumed to be equal to ${\mathbb P}^{M-2}$) is of
codimension $(i-1)$ in ${\mathbb P}^{M-2}$, but for some irreducible
component $B$ of this set we have $p_i|_B\equiv 0$,
and moreover if $i=M-2$, then $B$ is a curve of degree at least two.\vspace{0.1cm}

Obviously, Theorem 0.2 is implied by the following fact.\vspace{0.1cm}

{\bf Proposition 3.2.} {\it The following inequality holds:}
$$
\mathop{\rm codim}({\cal B}_i\subset{\cal P})\geqslant 2M.
$$

{\bf Proof.} Using the method that was applied in
\cite{Pukh98b} (see also \cite[Chapter 3, Subsection 1.3]{Pukh_book_13a}),
for  $i=1,\dots,k$ we obtain the estimate
$$
\mathop{\rm codim}({\cal B}_i\subset{\cal P})\geqslant\left(\begin{array}{c}
M+1-i\\
2\\\end{array}\right)
$$
(recall that $\delta(i)=2$ for $i=1,\dots,k$). The minimum of the right
hand sides is attained at $i=k$ and it is easy to check that
$$
\left(\begin{array}{c}M+1-i\\2\\\end{array}\right)-2M=
\frac12[(M-(2k+3))M+k^2-k]> 0.
$$
Therefore, we may assume that $i\geqslant k+1$, so that
$\delta(i)\geqslant 3$. Now we use the technique that was developed in
\cite{Pukh01} (see also \cite[Chapter 3, Section 3]{Pukh_book_13a}). Let
${\cal B}_{i,b}\subset{\cal P}$ be the set of tuples such that the
closed set (\ref{23.11.2016.1}) is of codimension $(i-1)$, and moreover,
there is an irreducible component $B$ of this set, such that
$$
\mathop{\rm codim}(\langle B\rangle\subset{\mathbb P}^{M-2})=b,
$$
$b\in\{0,1,\dots,i-1\}$,  $b\neq M-3$, and $p_i|_B\equiv 0$. Since
$$
{\cal B}_i=\bigcup^{i-1}_{b=0}{\cal B}_{i,b}
$$
(the condition $b\neq M-3$ for  $i=M-2$ is meant, but not shown, in order
for the formula not to be ugly), it is sufficient to show the inequality
$$
\mathop{\rm codim}({\cal B}_{i,b}\subset{\cal P})\geqslant 2M
$$
for $i\geqslant k+1$,  $b\in\{0,\dots,i-1\}$,  $b\neq M-3$.\vspace{0.1cm}

Now the technique of good sequences and associated subvarieties, which we do not
give here, see \cite{Pukh01} or
\cite[Chapter 3, Section 3]{Pukh_book_13a}, gives the estimate (taking into
account the dimension of the Grassmanian of linear subspaces of codimension $b$
in ${\mathbb P}^{M-2}$)
$$
\mathop{\rm codim}({\cal B}_{i,b}\subset{\cal P})\geqslant M(2b+3)-2b^2-6b-5.
$$
The right hand side of this inequality, considered as a function on the set
$\{0,\dots,i-1\}$, can decrease or increase or first increase and then decrease.
In any case the minimum of the right hand side is attained either at
$b=0$ (and equals $3M-5\geqslant 2M$), or at $b=i-1$ (if $i\leqslant M-3$) or
$b=M-4$ (if $i=M-2$), when it is also not smaller than $2M$.\vspace{0.1cm}

Q.E.D. for Theorem 0.2.

\newpage


\begin{flushleft}
Department of Mathematical Sciences,\\
The University of Liverpool
\end{flushleft}

\noindent{\it pukh@liverpool.ac.uk}

\end{document}